\documentclass[11pt, reqno]{article}
\usepackage{amsfonts}
\usepackage{amsmath}
\usepackage{amsthm}
\usepackage{amssymb}
\usepackage{lastpage}
\usepackage{fancyhdr}
\usepackage{fancybox}
\usepackage{mathrsfs}

\def\b{\mathbb }
\def\phi{\varphi }
\def\epsilon{\varepsilon}

\swapnumbers

\theoremstyle{plain}
\newtheorem{theorem}{Theorem}[section]
\newtheorem{corollary}[theorem]{Corollary}
\newtheorem{lemma}[theorem]{Lemma}
\newtheorem{proposition}[theorem]{Proposition}
\theoremstyle{definition}
\newtheorem{definition}[theorem]{Definition}
\newtheorem{remark}[theorem]{Remark}
\newtheorem{remarks}[theorem]{Remarks}

\numberwithin{equation}{section}

\begin{document}

\title{Limit theorems for radial random walks on $p\times q$-matrices as $p$ tends to infinity}
\author{Margit R\"osler\\
Institut f\"ur Mathematik, TU Clausthal\\
Erzstr. 1\\
D-38678 Clausthal-Zellerfeld, Germany\\
e-mail: roesler@math.tu-clausthal.de\\
and\\
Michael Voit\\
Fachbereich Mathematik, Universit\"at Dortmund\\
          Vogelpothsweg 87\\
          D-44221 Dortmund, Germany\\
e-mail: michael.voit@math.uni-dortmund.de}
\date{}
\maketitle

\begin{abstract}
The  radial probability measures on $\b R^p$ are in a one-to-one correspondence
with 
probability measures on $[0,\infty[$ by
 taking images of measures w.r.t.~the  Euclidean norm mapping. For fixed
 $\nu\in M^1([0,\infty[)$ and each dimension $p$, we consider i.i.d.~$\b
 R^p$-valued random variables $X_1^p,X_2^p,\ldots$ with radial laws
 corresponding to $\nu$ as above.
We derive weak and strong laws of large numbers as well as a
large deviation principle for the Euclidean length processes
$S_k^p:=\|X_1^p+\ldots+X_k^p\|$ as $k,p\to\infty$ in suitable ways.
\textbf{}
In fact, we  derive these results in a higher rank setting, where  
$\b R^p$ is replaced by the space of $p\times q$ matrices and  $[0,\infty[$
by the cone $\Pi_q$ of positive semidefinite matrices. Proofs are based on the
fact that the $(S_k^p)_{k\ge 0}$ form Markov chains on the cone
whose transition probabilities are given in terms Bessel functions $J_\mu$ of matrix argument
with an index  $\mu$ depending on $p$. The limit theorems  follow from new
 asymptotic results for the $J_\mu$  as $\mu\to \infty$. Similar results are also proven for
certain Dunkl-type Bessel functions.
\end{abstract}


\section{Introduction}

This paper has its origin in the  following problem:
Let $\nu\in M^1([0,\infty[)$ be a  probability measure. For each 
dimension $p\in\b N$ consider the time-homogeneous random walk
$(S_k^p)_{k\ge0}$ on $\b R^p$ which starts at time $k=0$ at $0\in\b R^p$ and makes  a random jump at each time step with  uniformly distributed
direction and a size  with distribution $\nu$ where the sizes and
directions are independent of each other and of the earlier ones. As the
distributions of the $ S_k^p$ are radial, we 
study  $\|S_k^p\|$ with
  the usual Euclidean norm $\|.\|$. Now let $(p_k)_{k\in\b N}\subset\b N$ be
a sequence of dimensions with $k\to\infty$. Our aim is to find limit
theorems for the $ [0,\infty[$-valued random variables $ \|S_k^{p_k}\|$ as
$k\to\infty$ for suitable sequences $(p_k)_{k\in\b N}\subset\b N$ of
dimensions  and suitable
measures $\nu$. This is an interesting question even for point
measures $\nu=\delta_r$ for $r>0$. We only have to exclude the trivial case
$\nu=\delta_0$, which we shall do from now on.

To get a first feeling for possible results, assume that $\nu$ has second
moment
$\sigma^2(\nu):=\int_0^\infty r^2\> d\nu(r)\in ]0,\infty[$. For each dimension
$p$ there is a unique radial measure $\nu_\rho\in M^1(\b R^p)$ with $\nu$ as
its radial part,
 i.e., for the norm mapping $\phi_p:\b
R^p\to [0,\infty[$, $x\mapsto \|x\|$, we have $\phi_p(\nu_p)=\nu$. We then may
realize the random walk $(S_k^p)_{k\ge0}$  as $S_k^p=\sum_{l=0}^k X_l^p$ for
i.i.d.~$\b R^d$-valued random variables $X_k^p$ with laws $\nu_p$ which admit
second moments.
The classical CLT on $\b R^p$ and the well-known relation between the standard
normal distribution on  $\b R^p$ and the $\chi^2$-distribution $\chi^2_p$ with
$p$ degrees of freedom imply after some short computation that for fixed $p$
and $k\to\infty$, the variables $\frac{p}{k\cdot \sigma^2(\nu)} \| S_k^p\|^2 $
tend in distribution to   $\chi^2_p$. Moreover, as $Z_p/p$ tends to 1 in probability for  $\chi^2_p$-distributed
random variables $Z_p$,  we obtain that 
$  \| S_k^p\|/\sqrt k\to \sqrt{ \sigma^2(\nu)} $
in probability if we first take $k\to\infty$ and then $p\to\infty$. 
We already observed in \cite{V1} that this result remains correct for other
combinations of $k,p\to\infty$:

\begin{theorem}\label{LLN1}
Assume that $\nu\in M^1([0,\infty[)$ has the second
moment
$\sigma^2(\nu)\in ]0,\infty[$.  Then for each
sequence  $(p_k)_{k\in\b N}\subset\b N$ of dimensions with $\lim_{k\to\infty}
p_k=\infty$, 
$$ \| S_k^{p_k}\|/\sqrt k\to \sqrt{ \sigma^2(\nu)} \quad\quad \text{in probability.} $$
\end{theorem}

One purpose of this paper is to prove an associated strong law. For simplicity
we  will assume that $\nu$ has a compact support.

\begin{theorem}\label{SLLN1}
Assume that $\nu\in M^1([0,\infty[)$ has compact support. Let  $(p_k)_{k\in\b
  N}\subset\b N$ and $(n_k)_{k\in\b N}\subset\b N$ sequences of dimensions and
time steps with the following properties:
\begin{enumerate}\itemsep=-1pt
\item[\rm{(1)}] $\lim_{k\to\infty}{p_k}/{k^a}=\infty$ for all $a\in\b N;$
\item[\rm{(2)}] $\lim_{k\to\infty}{p_k}/(n_k^2 (\ln k)^2)=\infty;$
\item[\rm{(3)}] $\lim_{k\to\infty}n_k/(\ln k)^2 =\infty$.
\end{enumerate}
  Then
$$ \| S_{n_k}^{p_k}\|/\sqrt{n_k}\to \sqrt{ \sigma^2(\nu)} \quad\quad\text{almost surely.} $$
\end{theorem}

For the case $n_k=k$,  only condition
(1) on the dimensions remains, i.e., the dimensions have to grow faster than any polynomial.
Unfortunately,  we are not able to get rid of this  strong growth
condition.
We shall discuss the conditions also in Section 4 below. Besides these laws of
large numbers we shall also derive a large deviation principle for
$S_{k}^{p_k}$ in Section 5 under the condition that $p_k$ grows faster than exponentially.

Theorems \ref{LLN1} and  \ref{SLLN1} and, in part, also this large deviation
principle will appear  as special
cases of extensions of these results in two directions. 

\bigskip 
The first extension concerns a higher rank setting.
We consider the following  
geometric situation: 
For fixed dimensions  $p,q\in \b N$ let 
$M_{p,q} = M_{p,q}(\b F)$  denote the space of $p\times q$-matrices over one of the division algebras
 $\b F = \b R, \b C$ or the quaternions $\b H$ with real dimension $d=1,2$ or $4$ respectively. This is a Euclidean vector space of (real) dimension $dpq$ with scalar product
$\langle x,y\rangle = \mathfrak R {tr}(x^*y)$ where $x^* := \overline x^t$, 
$\mathfrak R t := \frac{1}{2}(t+ \overline t)$ is the real part of $t\in \b F$,  and
${tr}$  is the trace in $M_{q}:= M_{q,q}.$
A  measure  on $M_{p,q}$ is called radial if it is 
invariant under the action of the unitary group $U_p= U_p(\b F)$ 
by left multiplication,
$U_p\times M_{p,q} \to M_{p,q}\,, \quad (u,x) \mapsto ux$.
This action is  orthogonal w.r.t.~the scalar product above, and, by uniqueness
of the polar decomposition, two matrices $x,y\in M_{p,q}$
belong to  the same $U_p$-orbit if and only if $x^*x = y^*y$. 
Thus the  space $M_{p,q}^{U_p}$ of $U_p$-orbits in $M_{p,q}$ is naturally parameterized by the 
cone $\Pi_q = \Pi_{q}(\b F)$ of positive semidefinite $q\times q$-matrices over
$\b F$. We
 identify $M_{p,q}^{U_p}$ with  $\Pi_q$ via $U_qx\simeq (x^*x)^{1/2}$, i.e., the
 canonical projection $M_{p,q}\to M_{p,q}^{U_p}$ will be realized as the mapping
$$ \phi_p: M_{p,q} \to \Pi_{q},\quad x\mapsto (x^*x)^{1/2}.$$ 
The square root is used here in order to ensure for  $q=1$ and  $\b F= \b R$ that 
 $\Pi_1=[0,\infty[$ and $\phi_p(x)=\|x\|$, i.e. the setting above appears.
By  taking images of measures,  the mapping $\phi_p$ induces a Banach space
isomorphism between the space $M_b^{U_q}(M_{p,q})$  of all bounded radial 
Borel measures on $M_{p,q}$ and  the space $M_b(\Pi_q)$ of bounded Borel measures on the cone $\Pi_q$.  In particular, for
each probability measure $\nu\in M^1(\Pi_q)$ there is a unique radial probability
measure $\nu_p\in M^1(M_{p,q})$ with $\phi_p(\nu_p)=\nu$. 
We shall say that $\nu\in M^1(\Pi_q)$ admits a second
moment if $\int_{\Pi_q}\|s\|^2 \> d\nu(s)<\infty$ where again, $\|s\|= (tr{s^2})^{1/2}$ is the Hilbert-Schmidt norm.  In this case, the second moment of $\nu$ is defined as the matrix-valued integral
$$\sigma^2(\nu):=\int_{\Pi_q} s^2 \> d\nu(s)\in \Pi_q.$$
With these notions, we shall derive the following generalizations of 
 Theorems \ref{LLN1} and \ref{SLLN1}:

\begin{theorem}\label{LLN2}
Let $\nu\in M^1(\Pi_q)$ be a probability measure with finite second moment
$\sigma^2(\nu)\in \Pi_q$. For each dimension $p\in\b N$ consider the unique
$U_p$-invariant probability measure $\nu_p\in M^1(M_{p,q})$ with
$\phi_p(\nu_p)=\nu$. Furthermore, let
$(X_l^p)_{l\in\b N}$ be a sequence of i.i.d.~$M_{p,q}$-valued random variables
with law $\nu_p$.
Then for each sequence $(p_k)_{k\in\b N}\subset \b N$ of dimensions with
$\lim_{k\to\infty} p_k=\infty$, 
$$\frac{1}{\sqrt k} \phi_{p_k}\bigl(\sum_{l=1}^k X_l^{p_k}\bigr) \to
\sqrt{\sigma^2(\nu)}\in\Pi_q
\quad\text{ in probability}.$$ 
\end{theorem}

\begin{theorem}\label{SLLN2}
Let $\nu\in M^1(\Pi_q)$ be a probability measure with compact support.
  For each dimension $p\in\b N$ consider the unique
$U_p$-invariant probability measure $\nu_p\in M^1(M_{p,q})$ with
$\phi_p(\nu_p)=\nu$. Furthermore, let
$(X_l^p)_{l\in\b N}$ be a sequence of i.i.d.~$M_{p,q}$-valued random variables
with law $\nu_p$.
Let  $(p_k)_{k\in\b
  N}\subset\b N$ and $(n_k)_{k\in\b N}\subset\b N$ sequences of dimensions and
time steps  with the following properties:
\begin{enumerate}\itemsep=-1pt
\item[\rm{(1)}] $\lim_{k\to\infty}{p_k}/{k^a}=\infty$ for all $a\in\b N;$
\item[\rm{(2)}] $\lim_{k\to\infty}{p_k}/(n_k^2 (\ln k)^2)=\infty;$
\item[\rm{(3)}] $\lim_{k\to\infty}n_k/(\ln k)^2 =\infty$.
\end{enumerate}
  Then
$\, \phi_{p_k}\bigl(\sum_{l=1}^{n_k} X_l^{p_k}\bigr)/{\sqrt{n_k}}$
tends  to $\sqrt{\sigma^2(\nu)}\in\Pi_q$ almost surely.
\end{theorem}

We next turn to a further generalization of these theorems.
Consider again the Banach space isomorphism between  $M_b^{U_q}(M_{p,q})$
and  $M_b(\Pi_q).$ The usual group convolution on $M_{p,q}$ induces 
a Banach-$*$-algebra-structure on $M_b(\Pi_q)$ such that this isomorphism becomes a
probability-preserving Banach-$*$-algebra isomorphism. 
 The space $\Pi_q$ together with this new convolution becomes a commutative orbit hypergroup; see \cite{J} and \cite{BH}
for a general background and \cite{R2}  for our specific
 example. It follows from  Eq.~(3.5) and Corollary 3.2 of \cite{R2} that in case $p\ge 2q$, the convolution product of two point measures on
$\Pi_q$ induced from $M_{p,q}$ is given by
\begin{equation}\label{def-convo-groupcase}
(\delta_r *_\mu \delta_s)(f) := 
\frac{1}{\kappa_{\mu}}\int_{D_q} f\bigl(\sqrt{r^2 + s^2 + svr + rv^*\!s}\,\bigr)\,
\Delta(I-vv^*)^{\mu-\rho}\, dv
\end{equation}
with $\mu:=pd/2$,
$\rho:=d\bigl(q-\frac{1}{2}\bigr) +1$,
$$ D_{q} := \{ v\in M_{q}: v^*v < I\}$$
(where $ v^*v <I$ means that $I- v^*v$ is strictly positive definite),
and with the normalization constant
\begin{equation}\label{kappa}
 \kappa_\mu := \int_{D_q} \Delta(I-v^*v)^{\mu-\rho} dv. 
\end{equation}
The convolution of arbitrary measures is just given by bilinear, weakly
continuous extension.

It was observed in  \cite{R2} that Eq.~(\ref{def-convo-groupcase})  defines a
commutative hypergroup actually for all indices $\mu\in\b R$ with $\mu>\rho-1$. In all cases, $0\in
\Pi_q$ is the  identity of the hypergroup and the involution is given by  the identity mapping. These hypergroup structures are closely related 
with a product formula for Bessel functions of index $\mu$  on the matrix cone $\Pi_q$ and 
are therefore called Bessel hypergroups on $\Pi_q$. Indeed, the hypergroup characters are given in terms of matrix
Bessel functions $J_\mu$. We refer to the monograph \cite{FK} for Bessel functions on cones,  and to \cite{R2} for
the particular details.
For general indices $\mu$, the Bessel hypergroups on $\Pi_q$ do not
have a nice geometric (orbit) interpretation as in the cases $\mu=pd/2$ with integral $p$, but nevertheless the notion of random walks on these
hypergroups is meaningful in the general cases just as well.

\begin{definition}
 Fix $\mu>\rho-1$ and a probability measure $\nu\in M^1(\Pi_q)$.
A Bessel random walk $(S_n^\mu)_{n\ge0}$ on $\Pi_q$ of index $\mu$ and with law $\nu$
 is a time-homogeneous Markov chain 
on $\Pi_q$   with  $S_0^\mu=0$ and   transition probability
$$P(S_{n+1}^\mu\in A| S_n^\mu=x) \> =\> (\delta_x*_\mu \nu)(A)$$
for $x\in\Pi_q$ and Borel sets $A\subset \Pi_q$. 
\end{definition}

This notion is quite common on hypergroups (see \cite{BH}) and was in particular
used  in \cite{V2} for  Bessel hypergroups on matrix cones and already in \cite{K} for the one-dimensional
case $q=1$ and $\b F=\b R$. The notion has its origin in the following well-known fact
for  the orbit cases $\mu=pd/2$ with $p\in\b N$: If we fix a radial measure $\nu_p\in
M^1(M_{p,q})$ and consider a sequence of i.i.d.~$M_{p,q}$-valued random variables
$(X_l^p)_{l\in\b N}$ 
with law $\nu_p$, then $ \bigl(\phi_{p}\bigl(\sum_{l=1}^k X_l^{p}\bigr)\bigr)_{k\ge0}$
is a random walk  on $\Pi_q$ of index $\mu$  with law $\phi_p(\nu_p)$.
Having this  in mind, we can state generalizations of Theorems \ref{LLN2} and
\ref{SLLN2} for such random walks on $\Pi_q$ for  indices $\mu\to \infty$
and time steps $k\to\infty$. This will be done in Section 4 where we state and
prove our results in this generality. The preceding limit results will then appear
just as special cases.

The proofs of the limit results in Section 4 are roughly as follows: As the characters of the
Bessel hypergroups on $\Pi_q$ can be expressed in terms of Bessel functions
$J_\mu$, the multidimensional Hankel transform on $\Pi_q$ is just the hypergroup
Fourier transform, and we can easily write down these transforms of the
distributions  of the $S_n^\mu$. On the other hand, we shall derive 
several uniform limit results for $J_\mu(\mu x)$ as $\mu\to\infty$. These results  imply  that the Hankel transforms tend to
 Laplace transforms of these distributions, which leads  to the stated
 limit theorems.
We point out that the direct proofs of Theorems \ref{LLN2} and \ref{SLLN2} are precisely
the same as in the slightly more general setting adopted in our paper.

The organization of this paper is as follows: In Section 2 we recapitulate some known
results about Bessel functions and Bessel convolutions on matrix cones from
\cite{FK}\cite{FT}\cite{H}, and \cite{R2}. The central part of the paper is Section 3, where  we 
present several uniform asymptotic results for $J_\mu(\mu x)$ as $\mu\to\infty$. Except for partial results proven by one of
the authors already in \cite{V1} for $q=1$, these results seem to be new even in the one-variable case $q=1$. This is surprising as in the classical monograph \cite{W} a complete
Chapter is devoted to  $J_\mu(\mu x)$ with $\mu\to\infty$. In Section 4, the asymptotic
results from Section 3 are transferred to certain classes of Dunkl-type Bessel functions associated with the root system $B_q$. Finally, the results of Section 3 
are used as a basis for the proofs of the laws of large numbers
in Section 5 and the large deviation principle in Section 6.

\section{Bessel functions and Bessel hypergroups on matrix cones}

In this section we collect some known facts about Bessel functions on matrix
cones and the associated Bessel hypergroups. The material is mainly taken from
 \cite{FK}  and  \cite{R2}. We also refer to the fundamental work \cite{H} of
 Herz, to \cite{Di}  and to \cite{FT}. 

\subsection{Bessel functions associated with matrix cones}

Let $\b F$ be one of the real division algebras $\b F= \b R, \b C$ 
or  $\b H$ with real dimension $d=1,2$ or $4$ respectively.
Denote the usual conjugation in $\b F$ by $t\mapsto \overline t$,
the real part  of $t\in \b F$ by $\mathfrak R t = \frac{1}{2}(t + \overline
t)$, and by $|t| = (t\overline t)^{1/2}$ its norm.

For $p,q\in \b N$ we denote  by $M_{p,q}:=M_{p,q}(\b F)$ 
the vector space of all  $p\times q$-matrices
 over $\b F$ and put $M_q:= M_q(\b F):=   M_{q,q}(\b F)$ for abbreviation. Let further
\[H_q = H_q(\b F)= \{x\in M_{q}(\b F): x= x^*\}\]
the space of Hermitian $q\times q$-matrices over $\b F$. All these spaces are real Euclidean
 vector spaces with scalar product 
 $\langle x, y\rangle := \mathfrak R {tr}(x^*y)$  and 
the associated norm $\|x\|=\langle x, x\rangle^{1/2}$. Here 
$x^* := \overline x^t$ and ${tr}$ denotes the trace.
The dimension of $H_q$ is given by  $dim_{\b R}H_q:= q + \frac{d}{2}q(q-1)$.
Let further
\[\Pi_q:=\{x^2:\> x\in H_q\} = \{ x^*x: x\in H_q\}\]
be the set of all positive semidefinite
matrices in $H_q$, and $\Omega_q$ its topological 
interior which consists of all strictly positive
definite matrices. $\Omega_q$ is a symmetric cone, i.e. an open convex
 cone which is self-dual and   whose  linear automorphism group acts transitively; see \cite{FK} for details.

To define the Bessel functions associated with the symmetric cone $\Omega_q$ we first  introduce their basic building blocks,
the so-called spherical polynomials. These are just 
 the polynomial spherical functions of $\Omega_q$ considered as a Riemannian symmetric space. They are indexed
by partitions $\lambda = (\lambda_1 \geq \lambda_2\geq \ldots \geq \lambda_q)
\in \b N_0^q$ (we write $\lambda \geq 0$ for short) and are given by 
\[ \Phi_\lambda (x) = \int_{U_q} \Delta_\lambda(uxu^{-1})du, \quad x\in H_q\]
where $du$ is the normalized Haar measure of $U_q$ and $\Delta_\lambda$ is the power function 
\[ \Delta_\lambda(x) := \Delta_1(x)^{\lambda_1-\lambda_2} \Delta_2(x)^{\lambda_2-\lambda_3} \cdot\ldots\cdot \Delta_q(x)^{\lambda_q}\quad (x\in H_q).\]
The $\Delta_i(x)$ are the principal minors of the determinant $\Delta(x)$, see 
\cite{FK} for details. 
There is a renormalization
$Z_\lambda = c_\lambda \Phi_\lambda$
with constants $c_\lambda >0$ depending on the underlying cone
such that
\begin{equation}\label{power-tr}
(tr \,x)^k \,=\, \sum_{|\lambda|=k} Z_\lambda(x)
\quad\quad\text{for}\quad
 k\ge0;
\end{equation}
see Section XI.5.~of \cite {FK} where these $Z_\lambda$ are called
zonal polynomials. By construction, the $Z_\lambda$ are invariant under
conjugation by $U_q$ and thus
 depend only on the eigenvalues of their argument.
More precisely, for $x\in H_q$ with eigenvalues $\xi = (\xi_1, \ldots, \xi_q)\in
\b R^q$, one has
\begin{equation}\label{identjack}Z_\lambda(x) = C_\lambda^\alpha(\xi) \quad \text{with}\quad \alpha = \frac{2}{d}\end{equation}
where the $C_\lambda^\alpha$ are the Jack polynomials of index
$\alpha$ in a suitable normalization (c.f. \cite {FK}, \cite{Ka}, \cite{R2}).
The Jack polynomials $C_\lambda^\alpha$ are homogeneous of degree $|\lambda|$ and symmetric in their arguments.

The matrix Bessel functions associated with the cone $\Omega_q$ are defined
as  $_0F_1$-hyper\-geometric series in terms of the $Z_\lambda$, namely
\begin{equation}\label{power-j}
 J_\mu(x) =
 \sum_{\lambda\geq 0} \frac{(-1)^{|\lambda|}}{(\mu)_\lambda|\lambda|!} Z_\lambda(x), 
\end{equation}
where for $\lambda=(\lambda_1, \ldots \lambda_q)\in \b N_0^q$,
 the generalized Pochhammer symbol $(\mu)_\lambda$ is given by
\[ (\mu)_\lambda =\, (\mu)_\lambda ^{2/d} \quad \text{where }\, (\mu)_\lambda^\alpha := \,\prod_{j=1}^q \bigl(\mu-\frac{1}{\alpha}(j-1)\bigr)_{\lambda_j} \quad (\alpha >0), \]
and $\mu\in \b C$ is an index satisfying $(\mu)_\lambda^\alpha \not= 0$ for all $\lambda \geq 0.$ If $q=1,$ then $\Pi_q=\b R_+$  and the Bessel function $\mathcal J_\mu$ is independent of $d$ with
\[ J_\mu\bigl(\frac{x^2}{4}\bigr) = j_{\mu-1}(x) \]
where $\, j_\kappa(z) = \,_0F_1(\kappa +1; -z^2/4)\,$
is the usual  modified Bessel function in one variable.

\subsection{Bessel hypergroups on matrix cones}

Hypergroups are convolution structures which generalize locally compact
groups insofar as the convolution product of two point measures is
in general not a point measure again, but just a probability measure on the underlying space. More precisely, a hypergroup $(X,*)$ is a locally compact Hausdorff space $X$ together with a convolution $*$
 on $M_b(X)$ (the regular bounded Borel measures on $X$),
such that $(M_b(X),*)$ becomes a Banach algebra, where  $*$ is
weakly continuous, probability preserving and  preserves compact
supports of measures.
 Moreover, one requires an identity $e\in X$  
with  $\delta_e*\delta_x=
\delta_x*\delta_e=\delta_x$ for $x\in X$, as well as
  a continuous involution $x\mapsto\bar x$ on $X$ such that for all $x, y\in X$,
$e\in supp(\delta_x*\delta_y)$ is equivalent to $ x=\bar y$,
and
$\delta_{\bar x} * \delta_{\bar y} =(\delta_y*\delta_x)^-$.
Here for $\mu\in M_b(X)$, the measure  $\mu^-$ is given by
$\mu^-(A)=\mu(A^-)$ for Borel sets $A\subset X$.
A hypergroup $(X,*)$ is called commutative if and only if so is the
convolution $*$. Thus for a commutative hypergroup $(X,*)$, the measure space
$M_b(X)$ becomes a commutative Banach-$*$-algebra with identity $\delta_e$.  Notice that due to its weak continuity, the convolution of measures on a hypergroup is uniquely determined by the convolution product of point measures. 

On a commutative hypergroup $(X,*)$ there exists a (up to a multiplicative factor) unique Haar measure $\omega$, i.e. $\omega$ is a positive Radon measure on $X$ satisfying 
\[ \int_X \delta_x*\delta_y(f)d\omega(y)  = \int_X f(y)d\omega(y)  \quad \text{for all } \, x\in X, \, f\in C_c(X).\]
The decisive object for harmonic analysis on a commutative hypergroup is its
dual space, which is defined by
\[ \widehat X :=\{\phi\in C_b(X): \phi\not= 0, \,\phi(\overline x) = \overline{\phi(x)}, \, \delta_x*\delta_y(\phi) = \phi(x)\phi(y) \, \text{ for all }\, x,y\in X\}.\] 
The elements of $\widehat X$ are also called characters. As in the case of LCA groups, the dual of a commutative hypergroup is a locally compact Hausdorff space with the topology of locally uniform 
convergence and can be identified with the symmetric spectrum of the convolution algebra $L^1(X,\omega).$

\medskip
The following theorem contains some of the main results of  \cite{R2}.

\begin{theorem}\label{main2}
Let $\mu\in\b R$ with $\mu > \rho-1.$ Then
\parskip=-1pt
\begin{enumerate}\itemsep=-1pt
\item[\rm{(a)}] The assignment
\begin{equation}\label{def-convo}
(\delta_r *_\mu \delta_s)(f) := 
\frac{1}{\kappa_\mu}\int_{D_q} f\bigl(\sqrt{r^2 + s^2 + svr + rv^*\!s}\,\bigr)\,
\Delta(I-vv^*)^{\mu-\rho}\, dv, \quad f\in C(\Pi_q)
\end{equation}
with  $\kappa_\mu$ as in (\ref{kappa}),  defines a commutative hypergroup structure on $\Pi_q$ with 
neutral element $0\in \Pi_q$ and the identity mapping as involution.
The support of $\delta_r*_\mu\delta_s$ satisfies
\[\text{supp}(\delta_r*_\mu\delta_s) \subseteq  
\{t\in \Pi_q: \|t\|\leq \|r\|+\|s\|\}.\]
\item[\rm{(b)}] A Haar measure of the hypergroup 
$\Pi_{q,\mu}:= (\Pi_q, *_\mu)$ is given by
\[ \omega_\mu(f) = 
\frac{\pi^{q\mu}}{\Gamma_{\Omega_q}(\mu)}
\int_{\Omega_q} f(\sqrt{r}) \Delta(r)^\gamma dr\]
with $\,\gamma = \mu-\frac{d}{2}(q-1)-1.$
\item[\rm{(c)}]
The dual space of $\Pi_{q,\mu}$ is given by
\[\widehat{\Pi_{q,\mu}} = \,\{\phi_s: s\in \Pi_q\}\]
with 
$$\phi_s(r):= \mathcal J_\mu(\frac{1}{4}rs^2r)=\phi_r(s).$$
 The hypergroup $\Pi_{q,\mu}$ is self-dual via the homeomorphism $s\mapsto \phi_s$.
Under this identification of $\widehat{\Pi_{q,\mu}}$ with $\Pi_{q,\mu}\,$,
  the Plancherel measure on $\Pi_{q,\mu}\,$ is $(2\pi)^{-2\mu q}\omega_\mu$.
\end{enumerate}
\end{theorem}

\section{Estimates for Bessel functions of large indices}

We first recapitulate the following well known one-dimensional 
inequalities for the exponential function 
(see, for instance,  Sections 3.6.2 and 3.6.3 of \cite{Mi}):
\begin{equation}\label{mit0}
(1-z/r)^r \le  e^{-z} \quad\text{for}\quad r>0, z\in\b R,
\end{equation}
\begin{equation}\label{mit}
0\le e^{-z}-(1-z/r)^r \le z^2 e^{-z}/r \quad\text{for}\quad
r,z\in\b R,\> r\ge 1, \> |z|\le r;
\end{equation}
\begin{equation}\label{mit1}
(1+z/r)^r\le e^z \le (1+z/r)^{r+z/2}
\quad\text{for}\quad
r,z>0.
\end{equation}
These results have the following matrix-valued extension:

\begin{lemma}\label{lemma1}
\begin{enumerate}\itemsep=-1pt
\item[\rm{(1)}] 
For all $\mu>1$ and $v\in\sqrt \mu\cdot D_q\subset M_{q}$,
$$0\le e^{-\langle v,v\rangle}-\Delta(I-\frac{1}{\mu}v^*v)^\mu \>\le \>
\frac{1}{\mu} tr((v^*v)^2)\cdot e^{-\langle v,v \rangle}.$$
\item[\rm{(2)}] For all $\mu>0$ and $v\in M_{q}$,
$$\Delta(I+\frac{1}{\mu}v^*v)^\mu \le  e^{\langle v,v\rangle}\le
\Delta(I+\frac{1}{\mu}v^*v)^{\mu+m/2}$$
where $m \ge0$ is the maximal eigenvalue of $v^*v$.
\end{enumerate}
 \end{lemma}

\begin{proof}
\begin{enumerate}\itemsep=-1pt
\item[\rm{(1)}] 
The positive semidefinite matrix $v^*v$ may be written as
$$v^*v= u\cdot diag(a_1,a_2,\ldots,,a_q)\cdot u^*$$
with some $u\in U_q$ and the eigenvalues  $a_1,\ldots,a_q\in [0,m]$ 
of $v^*v$. Then
$$\Delta(I-\frac{1}{\mu}v^*v)^\mu=\prod_{k=1}^q(1-a_k/\mu)^\mu$$
and  $\langle v,v\rangle=tr(v^*v)= a_1+\ldots+a_q$. Using (\ref{mit}) and a
telescope sum argument, we obtain
\begin{align}
0 &\le  e^{-\langle v,v\rangle}-\Delta(I-\frac{1}{\mu}v^*v)^\mu
\notag\\
 &= \sum_{l=1}^q \Bigl[\prod_{k=1}^{l-1} (1-a_k/\mu)^\mu \cdot 
\bigl(e^{-a_l}-(1-a_l/\mu)^\mu\bigr)\cdot  \prod_{k=l+1}^{q}e^{-a_k}
\bigr]
\notag\\
 &\le \frac{1}{\mu}\sum_{l=1}^q\Bigl[   \prod_{k=1}^qe^{-a_k}\bigr] a_l^2
\>=\> \frac{1}{\mu} tr((v^*v)^2)\cdot e^{-\langle v,v \rangle}
\notag
\end{align}
as claimed.
\item[\rm{(2)}] is proven in the same way by use of (\ref{mit1}).
\end{enumerate} 
\end{proof}

Our first estimate for the Bessel functions 
$J_\mu$ as $\mu\to\infty$ will be based on the following  integral representation
of  $J_\mu$ for $\mu >\rho -1$ (see Eq.~(3.12) of \cite{R2}):
\begin{equation}\label{int-rep}
J_\mu(x^*x)=
\,\frac{1}{\kappa_\mu} \int_{D_q} e^{-2i\langle v, x\rangle} \Delta(I-v^*v)^{\mu-\rho}dv.
\end{equation}

\begin{proposition}\label{proposition1}
There exists a constant $C=C(q,d)>0$ such that for 
$\mu>\rho-1$ and all
$x\in M_{q}$,
$$\bigl|J_\mu(\mu x^*x) - e^{-\langle x,x \rangle}\bigr|\le C/\mu$$
and
$$| \mu^{dq^2/2} \kappa_\mu- \pi^{dq^2/2}|\le  C/\mu.$$
 \end{proposition}

\begin{proof} In a first step we obtain for $x\in M_q$,
\begin{align}
D:= &\,\Bigl|\int_{\sqrt\mu\cdot D_q}  e^{-i\langle v, x\rangle} 
\Delta(I-\frac{1}{\mu}v^*v)^{\mu-\rho}\> dv - \int_{M_{q}}  e^{-i\langle v,
  x\rangle}
 e^{-\langle v, v\rangle} \> dv\Bigr|
\notag\\
 &\le \int_{M_{q}\setminus (\sqrt\mu\cdot D_q)}e^{-\langle v, v\rangle}
 \> dv \>+\>
\int_{\sqrt\mu\cdot D_q} \Bigl(e^{-\langle v, v\rangle(1-\rho/\mu)}  - 
\Delta\bigl(I-\frac{1-\rho/\mu}{\mu-\rho}v^*v\bigr)^{\mu-\rho}\Bigr) \> dv
\notag\\& + \int_{\sqrt\mu\cdot D_q}e^{-\langle v, v\rangle}
\bigl(e^{\langle v, v\rangle \rho/\mu}-1\bigr) \> dv.
\notag
\end{align}
By
 Lemma \ref{lemma1}(1) for $\mu-\rho$ instead of $\mu$ and with the elementary
estimate
$$e^z-1\le (e^\rho-1)z \quad\text{for}\quad z\in [0,\rho]$$ we further obtain that for $\mu\geq 2\rho$,
\begin{align}\label{large-est} D\,\leq &\,\frac{C_1}{\mu} + \frac{1}{\mu-\rho}\int_{\sqrt\mu\cdot D_q}e^{-\langle
   v, v\rangle/2} tr((v^*v)^2) \> dv\, +\, \frac{(e^\rho-1)\rho}{\mu}
\int_{\sqrt\mu\cdot D_q}e^{-\langle
   v, v\rangle}\langle
   v,v\rangle \> dv\notag\\
 \leq & \, \frac{C_2}{\mu}
\end{align}
with suitable constants $C_1,C_2>0$.
We next observe that 
$M_{q}\simeq\b R^{dq^2}$ implies
\begin{equation}\label{int-rep-exp}
\int_{M_{q}}e^{-i\langle v, x\rangle} e^{-\langle v, v\rangle} \>
dv
\>=\>
\pi^{dq^2/2}\cdot e^{-\langle x, x\rangle/4}.
\end{equation}
Moreover, 
replacing $x$ by $(\sqrt\mu/2)\cdot x$ and $v$ by $(1/\sqrt\mu)\cdot v$ in
integral representation (\ref{int-rep}), we obtain
\begin{equation}\label{int-rep2}
J_\mu(\frac{\mu}{4}x^*x)=
\,\frac{1}{\mu^{dq^2/2} \kappa_\mu} \int_{\sqrt\mu\cdot D_q}
 e^{-i\langle v, x\rangle} 
\Delta(I-\frac{1}{\mu}v^*v)^{\mu-\rho}dv.
\end{equation}
We now conclude from (\ref{large-est}), (\ref{int-rep-exp}) and
(\ref{int-rep2}) that for  $\mu\ge2\rho$ and  $x\in M_{q}$,
$$\Bigl|\mu^{dq^2/2}\kappa_\mu\cdot J_\mu(\frac{\mu}{4}x^*x) -
\pi^{dq^2/2}\cdot e^{-\langle x, x\rangle/4}\Bigr|\le C_2/\mu.$$
 For $x=0$ we in particular
observe that
\begin{equation}\label{kappa-absch}
| \mu^{dq^2/2} \kappa_\mu- \pi^{dq^2/2}|\le  C_2/\mu.
\end{equation}
As $|J_\mu(\frac{\mu}{4}x^*x)|\le1$  by  (\ref{int-rep}),  it follows for  $\mu\ge2\rho$ and $x\in M_{q}$ that
\begin{align}
\bigl|&J_\mu(\frac{\mu}{4}x^*x) - e^{-\langle x,x \rangle/4}\bigr|
\notag\\ &\le
|J_\mu(\frac{\mu}{4}x^*x)|\cdot \bigl| 1- \frac{\mu^{dq^2/2}\cdot \kappa_\mu}{ \pi^{dq^2/2}}\bigr|
+\frac{1}{ \pi^{dq^2/2}}\bigl| \mu^{dq^2/2}\cdot \kappa_\mu J_\mu(\frac{\mu}{4}x^*x) - 
 \pi^{dq^2/2}e^{-\langle x, x\rangle/4}\bigr|
\notag\\ &\le
C_3/\mu
\notag
\end{align}
 with some constant $C_3>0$. Together with (\ref{kappa-absch}), this implies the statements of the proposition in case
 $\mu\ge2\rho$. Within the range $\rho -1 < \mu \leq 2\rho$, the proposition is  immediate in view of the estimate $|J_\mu(\mu x^*x)|\le1$ for all $x\in M_q$.
\end{proof}

In the following, we shall derive a variant of Proposition \ref{proposition1} which is based on the power
series (\ref{power-j}) and provides a good estimate for small arguments. 
We start with some basic inequalities for the zonal polynomials $Z_\lambda$:

\begin{lemma}\label{lemma20}
For all  partitions $\lambda\ge0$ and 
and  $y\in \Pi_q$,
$$|Z_\lambda(-y)|\le Z_\lambda(y).$$
\end{lemma}

\begin{proof}
We use the relation between the $Z_\lambda$ and  the Jack polynomials
$C_\lambda^\alpha$   in Section 2 and the well-known fact that the $C_\lambda^\alpha$ are
nonnegative linear combinations of monomials, see \cite{KS}.  This yields for the eigenvalues 
$\xi=(\xi_1,\ldots ,\xi_q)$ of $-y$ and $|\xi|:=(|\xi_1|,\ldots ,|\xi_q)|)$ of $y$ that
 $$|Z_\lambda(-y)|= |C_\lambda^\alpha(\xi)|\le
 C_\lambda^\alpha(|\xi|)=Z_\lambda(y).$$
\end{proof}

\begin{lemma}\label{lemma2}
For all  partitions $\lambda\ge0, \, \mu>\rho-1$ and $(\mu)_\lambda= (\mu)_\lambda^{2/d}$,
$$\Bigl|1-\frac{\mu^{|\lambda|}}{(\mu)_\lambda}\Bigr|\le 
dq \cdot 2^{dq(q-1)/2}\cdot\frac{ |\lambda|^2}{\mu}.$$
\end{lemma}

\begin{proof}
Consider $\, (\mu)_\lambda \>=\> \prod_{j=1}^q(\mu-\frac{d}{2}(j-1))_{\lambda_j}.$
In this product, each factor can be estimated below by $\mu-\frac{d}{2}(q-1)$. Moreover, precisely
\[ \bigl(0 + 1 + \ldots  + (q-1) \bigr)\lceil d/2\rceil \, = \, \frac{(q-1)q}{2}\cdot\lceil d/2 \rceil\,=:r\]
of these factors are smaller than $\mu$. As $\mu>\rho-1=d(q-1/2),$ this implies 
\begin{align}
(\mu)_\lambda &
\>\ge\> (\mu-\frac{d}{2}(q-1))^r \cdot \mu^{|\lambda|-r}
\notag\\ 
&\>\ge\> \bigl(\mu/2\bigr)^r  \mu^{|\lambda|-r} 
\>\ge\> 2^{-dq(q-1)/2} \cdot \mu^{|\lambda|},
\notag
\end{align}
and thus
\begin{equation}\label{absch-poch}
 \mu^{|\lambda|}/ (\mu)_\lambda \>\le\>  2^{dq(q-1)/2}.
\end{equation}

We now prove by
 induction on the length $k:=|\lambda|$ that
 for $\mu>\rho-1=d(q-1/2)$,
\begin{equation}\label{absch}
\Bigl|1-\frac{\mu^{|\lambda|}}{(\mu)_\lambda}\Bigr|\le 
\frac{dq}{2(\mu-d(q-1)/2)}\cdot 2^{dq(q-1)/2} |\lambda|^2
\end{equation}
which immediately implies the lemma. 
 In fact, for $k=0,1$, the left
hand side of (\ref{absch}) is equal to zero, while the right-hand side is nonnegative. 

For the induction step, consider a partition $\lambda$ of length $k\ge 2$. 
Then there is  a partition $\tilde \lambda$ with $|\tilde \lambda|=k-1$ for
which there exists precisely one $j=1,\ldots,q$ with
$\lambda_j=\tilde\lambda_j +1$ while all the other components are equal. Hence, if
we assume the inequality to hold for $ \tilde \lambda$ and use (\ref{absch-poch}), we obtain
\begin{align}
\Bigl|1&-\frac{\mu^k}{(\mu)_\lambda}\Bigr|=
 |1-\frac{\mu^{k-1}}{(\mu)_{\tilde\lambda}}
+\frac{\mu^{k-1}}{(\mu)_{\tilde\lambda}}
-\frac{\mu^k}{(\mu)_\lambda}
\Bigr|
\notag\\ &\le
\frac{dq}{\mu-d(q-1)/2}\cdot 2^{dq(q-1)/2-1} \cdot (k-1)^2
 +
\frac{\mu^{k-1}}{(\mu)_{\tilde\lambda}}\cdot 
\Bigl|1-\frac{\mu}{(\mu-d(j-1)/2)+\lambda_j-1}\Bigr|
\notag\\ &\le
\frac{dq}{\mu-d(q-1)/2}\cdot 2^{dq(q-1)/2-1} \cdot (k-1)^2
 +
 2^{dq(q-1)/2}
\cdot 
\Bigl|\frac{ -d(j-1)/2 +\lambda_j-1}{(\mu-d(j-1)/2)+\lambda_j-1}\Bigr|
\notag\\ &
\le
\frac{2^{dq(q-1)/2-1}}{\mu-d(q-1)/2}\cdot\Bigl( dq(k-1)^2 +dq+2k-2\Bigr)
 \notag\\ &
\le  
\frac{2^{dq(q-1)/2-1}}{\mu-d(q-1)/2}\cdot dq k^2
\notag
\end{align}
for $k\ge2$. This completes the proof.
\end{proof}

\begin{proposition}\label{proposition2}
There exists a constant $C=C(q,d)>0$ such that for
$\mu>2\rho$ and 
$y\in \Pi_q$,
$$\bigl|J_\mu(\mu y) - e^{-tr\, y }\bigr|\le
C\frac{(tr\, y)^2}{\mu}.$$
 \end{proposition}

\begin{proof}
Using the power series (\ref{power-j}) as well as (\ref{power-tr}) in terms of the
homogeneous polynomials $Z_\lambda$, we obtain
$$J_\mu(\mu y)-e^{-tr\, y}  = \sum_{\lambda\geq 0} 
\frac{1}{|\lambda|!}\Bigl( \frac{\mu^{|\lambda|}}{(\mu)_\lambda}-1\Bigr)
\cdot Z_\lambda(-y).$$
As
$$(\mu)_{(1,0,\ldots,0)}=\mu,\quad (\mu)_{(2,0,\ldots,0)}=\mu(\mu+1),
\quad (\mu)_{(1,1,0,\ldots,0)}=\mu(\mu-d/2), $$
 we may write this expansion as
$$J_\mu(\mu y)-e^{-tr\, y}  = R_2+R_3$$
with
$$R_2= \frac{1}{2}\Bigl(\Bigl(\frac{\mu^2}{\mu(\mu+1)}-1\Bigr)
Z_{(2,0,\ldots,0)}(-y)+\Bigl(\frac{\mu^2}{\mu(\mu-d/2)}-1\Bigr)
Z_{(1,1,0,\ldots,0)}(-y)\Bigr)$$
and
$$R_3=\sum_{k\geq 3} \frac{1}{k!} \sum_{|\lambda|= k} 
\Bigl( \frac{\mu^{k}}{(\mu)_\lambda}-1\Bigr)
\cdot Z_\lambda(-y).$$
Recall from  Lemma \ref{lemma20} that $|Z_\lambda(-y)|\le Z_\lambda(y)$ and $Z_\lambda(y) \geq 0$. 
Hence Eq.~(\ref{power-tr}) implies  for $|\lambda|=2$ that 
 $|Z_\lambda(-y)|\le (tr\, y)^2$. Therefore,
$|R_2|\le M_1 \frac{(tr\,y)^2}{\mu}$ with a suitable constant $M_1>0$.
Moreover, Lemmata \ref{lemma20} and  \ref{lemma2} imply that 
\begin{align}
|R_3 |  
&\> \le  \>\sum_{k\geq 3} \frac{1}{k!}\sum_{|\lambda|= k} M_2 \frac{k^2}{\mu} Z_\lambda(y)
\> =  \> \frac{M_2}{\mu} \sum_{k\geq 3} \frac{k^2}{k!} (tr\, y)^k
 \notag\\
 &\> \le  \> \frac{2M_2}{\mu} (tr\, y)^2 \sum_{k\geq 1} \frac{1}{k!} (tr\, y)^k
\> \le  \> \frac{2M_2}{\mu} (tr\, y)^2 e^{tr\, y}
\end{align}
with a constant $M_2>0$. In summary we have 
$$ |J_\mu(\mu y)-e^{-tr\, y} | \le M_3 \cdot\frac{(tr\, y)^2}{\mu}\cdot(1+e^{tr\, y}).$$
Together with the estimate of  Proposition \ref{proposition1} for large $y$, this yields the stated result. 
\end{proof}

Summarizing Propositions \ref{proposition1} and \ref{proposition2}, we obtain:

\begin{theorem}\label{theorem1}
There exists a constant $C=C(q,d)>0$ such that for
$\mu>2\rho$ and 
$y\in \Pi_q$,
$$\bigl|J_\mu(\mu y) - e^{-tr \, y }\bigr|\le
\frac{C}{\mu}\cdot  min(1, (tr\, y)^2).$$
\end{theorem}

We next turn to an estimate for $J_\mu(-\mu y)$ with $y\in \Pi_q$. In order to
simplify formulas, we replace the factor $\mu$ in the argument by $\mu-\rho$.

\begin{proposition}\label{proposition3}
There exists a constant $C=C(q,d)>0$ such that for
$\mu>2\rho$ and all
$x\in M_{q}$,
$$ e^{\langle x,x\rangle} \Bigl(1-\frac{C}{\mu}\|x\|^4 - H(x,\sqrt{\mu-\rho})\Bigr)
  \le J_\mu(-(\mu-\rho)x^*x)  \le e^{\langle x,x\rangle} (1+C/\mu)$$
where
$$H(x,r):=\int_{M_{q}\setminus r D_q} e^{-\|v-x\|^2}\> dv
\quad\quad{\rm for}\quad r>0.$$
\end{proposition}

\begin{proof}
We first conclude (by analytic continuation) from  integral representation (\ref{int-rep}) that
\begin{equation}\label{int-rep3}
J_\mu(-(\mu-\rho)x^*x)=
\,\frac{1}{(\mu-\rho)^{dq^2/2} \kappa_\mu} \int_{\sqrt{\mu-\rho}\cdot D_q}
 e^{2\langle v,x\rangle} 
\Delta(I-\frac{1}{\mu-\rho}v^*v)^{\mu-\rho}dv.
\end{equation}
Moreover, Proposition \ref{proposition1} implies that
\begin{equation}\label{eq2}
\Bigl|  \frac{1}{(\mu-\rho)^{dq^2/2} \kappa_\mu}  -
\pi^{-dq^2/2}\Bigr|=O(1/\mu).
\end{equation}
We next estimate the integral in (\ref{int-rep3}). For this we use Lemma
\ref{lemma1}(1) and observe that
\begin{align} 
\int_{\sqrt{\mu-\rho}\cdot D_q}
 e^{2\langle v, x\rangle} 
\Delta(I-\frac{1}{\mu-\rho}v^*v)^{\mu-\rho}\> dv &\le 
\int_{M_{q}} e^{2\langle v, x\rangle} e^{-\langle v, v \rangle} \> dv
\notag\\
&= \pi^{dq^2/2} e^{\langle x, x \rangle}. 
\notag
\end{align}
Together with Eq.~(\ref{int-rep3})    and   (\ref{eq2}) this yields
$$J_\mu(-(\mu-\rho)x^*x)\le (\pi^{-dq^2/2}+ O(1/\mu))\cdot
\pi^{dq^2/2}e^{\langle x, x \rangle}.$$
This proves the upper estimate as claimed.

For the lower estimate, we use Lemma \ref{lemma1}(1) again. We obtain
\begin{align} 
\int_{\sqrt{\mu-\rho}\cdot D_q}
 e^{2\langle v, x\rangle} 
\Delta(I-\frac{1}{\mu-\rho}v^*v)^{\mu-\rho}\> dv &\ge 
\int_{\sqrt{\mu-\rho}\cdot D_q } e^{2\langle v, x\rangle-\langle v, v
  \rangle}\cdot
\bigl( 1-\frac{1}{\mu-\rho}tr((v^*v)^2)\bigr) \> dv
\notag\\
&= \int_{M_{q}} e^{2\langle v, x\rangle-\langle v, v
  \rangle} \> dv  - I_1 - \frac{1}{\mu-\rho}I_2
\notag\\
&= \pi^{dq^2/2} e^{\langle x, x \rangle} - I_1 - \frac{1}{\mu-\rho}I_2
\notag
\end{align}
with
$$I_1 := \int_{M_{q}\setminus \sqrt{\mu-\rho}\cdot D_q} 
e^{2\langle v, x\rangle-\langle v, v \rangle}
\> dv$$
and
$$I_2 := \int_{ \sqrt{\mu-\rho}\cdot D_q} e^{2\langle v, x\rangle-\langle
  v, v \rangle} \cdot tr((v^*v)^2)
\> dv.$$
We have
$$I_1 = e^{\langle x, x \rangle} \cdot \int_{M_{q}\setminus \sqrt{\mu-\rho}\cdot D_q} 
 e^{-\|v-x\|^2}\> dv = 
e^{\langle x, x \rangle} \cdot H(x,\sqrt{\mu-\rho})$$
and 
$$I_2= e^{\langle x, x \rangle} \cdot \int_{\sqrt{\mu-\rho}\cdot D_q} 
e^{-\|v-x\|^2}  \cdot tr((v^*v)^2)
\> dv = e^{\langle x, x \rangle} \cdot O(\|x\|^4),$$
which finally leads to the lower estimate.
\end{proof}

\section{Estimates for Dunkl-type Bessel functions associated with root systems of type B}

There is a close connection between Bessel convolutions on the cone $\Pi_q$ 
and the theory of Dunkl operators associated with the root system  $B_q$ which is explained in  \cite{R2}. In this short section, we shall 
recall this connection and use it to obtain asymptotic relations between
certain classes of Dunkl-type Bessel functions, which can be expressed as generalized hypergeometric functions in terms of Jack polynomials. This section is independent of the remaining parts of this paper and may be skipped by readers interested in the probabilistic results only. Also, we shall not go into details of Dunkl theory, but refer the  reader to \cite{DX}, 
\cite{R1} and \cite{R3}.  For multivariable hypergeometric functions, see e.g. \cite{GR} and \cite{Ka}.  In the following, we always assume that $q\geq 2$.  For a reduced root system $R\subset \b R^q$ and a multiplicity function $k:R\to \b C$  (i.e. $k$ is  invariant under the action of the corresponding reflection group), 
we denote by $J_k = J_k^R$ the Dunkl-type Bessel function associated with $R$ and $k$. It is obtained from the Dunkl kernel by symmetrization with respect to the underlying reflection group. Dunkl-type Bessel functions generalize the spherical functions of Euclidean type symmetric spaces, which occur for crystallographic root systems and specific discrete values of $k$.  For the root system $A_{q-1}= \{\pm(e_i-e_j):\, i<j\}\subset \b R^q$, the multiplicity $k$ is 
a single complex parameter and if $k>0$, then the associated Dunkl-type Bessel function can be expressed as a generalized $_0F_0$-hypergeometric function,
\[ J_k^A(\xi,\eta) =\, _0F_0^\alpha (\xi,\eta) := \sum_{\lambda \geq 0} \frac{1}{|\lambda|!}\cdot 
\frac{C_\lambda^\alpha(\xi) C_\lambda^\alpha(\eta)}{C_\lambda^\alpha({\bf 1})}\quad \text{with }\, {\bf 1} = (1,\dots, 1),\, \alpha = 1/k, \]
due to relations (3.22) and (3.37)  of \cite{BF}. 
For the root system   $B_q = \{\pm e_i , \,\pm e_i \pm e_j: \, i <j \}$, the multiplicity is of the form $k=(k_1, k_2)$ where $k_1$ and $k_2$ are the values on the roots $\pm e_i$ and  $\pm e_i \pm e_j$ respectively. The associated Dunkl-type Bessel function is given by 
\[ J_k^B(\xi,\eta) = \, _0F_1^\alpha\bigl(\mu; \frac{\xi^2}{2}, \frac{\eta^2}{2}\bigr)
\quad \text{with }\, \alpha = 1/k_2, \, \mu = k_1 +(q-1)k_2 +1/2\]
where $\xi^2= (\xi_1^2, \ldots, \xi_q^2)$ and 
\[_0F_1^\alpha(\mu; \xi, \eta) := \, \sum_{\lambda\geq 0} \frac{1}{(\mu)_\lambda^\alpha |\lambda|!}\cdot \frac{C_\lambda^\alpha(\xi) C_\lambda^\alpha(\eta)}{C_\lambda^\alpha({\bf 1})}.\] 

Recall now that the characters of the hypergroup $\Pi_{q,\mu}$ on the matrix cone $\Pi_q$ are given by 
$\phi_s(r) = \mathcal J_\mu\bigl(\frac{1}{4}sr^2s\bigr).$ The conjugation action $x\mapsto uxu^{-1}$ 
of the unitary group $U_q=U_q(\b F)$ on $\Pi_q$ induces a new commutative hypergroup structure on the set
of possible eigenvalues of matrices from $\Pi_q$ ordered by size, i.e. the $B_q$-Weyl chamber
\[ \Xi_q = \{ \xi = (\xi_1,\ldots ,\xi_q) \in \b R^q: \xi_1\geq \ldots \geq \xi_q\geq 0\}. \]
This hypergroup on $\Xi_q$  depends on $d$ and $\mu$. Its characters 
are given by the functions 
\begin{equation}\label{Dunklchar} \psi_\eta(\xi)=\,\int_{U_q} \mathcal J_\mu\bigl(\frac{1}{4}\eta u\xi^2 u^{-1}\eta\bigr) du \,= \, J_{k(\mu,d)}^B(\xi, i\eta), \quad \eta\in \Xi_q\end{equation}
where $\,k(\mu, d) = \bigl(\mu - (d(q-1)+1)/2, \,d/2\bigr)$ (and elements from $\Xi_q$ are identified with diagonal matrices in the natural way). 
For details, see Section 4 of \cite{R2}.
The estimates for the  matrix Bessel functions $\mathcal J_\mu$ according to Theorem \ref{theorem1} 
imply the following estimate for the Dunkl-type Bessel function $J_{k(\mu,d)}^B$ as $\mu\to\infty$. 

\begin{corollary} There exists a constant $C=C(q,d)>0$ such that for $\mu >2\rho$ and $\xi, \eta\in \Xi_q$,
\[ \big\vert J_{k(\mu,d)}^B(2\sqrt{\mu}\,\xi, i\eta) \,-\, J_{d/2}^A(-\xi^2,\eta^2) \big\vert\,
  \leq\, \frac{C}{\mu} \cdot\min \bigl(1, (|\xi^2||\eta^2|)^2\bigr)\]
where $|\zeta|= (\zeta_1^2 + \ldots \zeta_q^2)^{1/2}$ denotes the standard Euclidean norm in $\b R^q$. 
\end{corollary}

\begin{proof}
For $k=k(\mu,d)$ we obtain by Eq.~(\ref{Dunklchar}) and 
Theorem \ref{theorem1} the estimate
\[\Big\vert J_k^B(2\sqrt{\mu}\, \xi,i\eta) \,-\, \int_{U_q} e^{-tr(\eta u \xi^2u^{-1}\eta)}du \Big\vert \,\leq \, \frac{C}{\mu} \cdot\min \bigl(1, S(\xi,\eta) \bigr)\]
where
\[ S(\xi,\eta) \,=\, \int_{U_q}\big\vert tr (\eta u\xi^2 u^{-1} \eta)\big\vert ^2 du \,= \, \int_{U_q}
\vert \langle \eta^2, u \xi^2 u^{-1}\rangle\vert^2 du\, \leq \, |\xi^2|^2|\eta^2|^2\]
by the Cauchy-Schwarz inequality. 
The spherical polynomials $Z_\lambda$  satisfy the product formula 
\[ \frac{Z_\lambda(r)Z_\lambda(s)}{Z_\lambda(I)}  = \int_{U_q} Z_\lambda(\sqrt{r}usu^{-1}\sqrt{r})du \quad \text{for }\,r,s\in \Pi_q\,,\]
see Prop. 5.5. of \cite{GR}. Thus by Eq.~(\ref{power-tr}) and 
(\ref{identjack}) we further obtain, with $\alpha = 2/d$, 
\begin{align}\label{harish} \int_{U_q} e^{-tr(\eta u \xi^2u^{-1}\eta)}du\,= &\, \sum_{\lambda\geq 0} 
\frac{1}{|\lambda|!} \int_{U_q} Z_\lambda (-\eta u \xi^2u^{-1}\eta)du \notag\\
=&\, \sum_{\lambda\geq 0} \frac{1}{|\lambda|!}\frac{Z_\lambda(-\xi^2)Z_\lambda(\eta^2)}{Z_\lambda(I)}
\,=\, _0F_0^\alpha (-\xi^2,\eta^2),\end{align}
which implies the assertion. 
\end{proof}

\begin{remarks}
\begin{enumerate}\itemsep=-1pt
\item[\rm{(1)}]  It is conjectured that the statement of this corollary  remains valid for arbitrary $d\in \b R$ with $d>0$. 
\item[\rm{(2)}] The integral on the left side of formula (\ref{harish}) is of Harish-Chandra type. If $\b F = \b C$, then by Theorem II. 5.35 of \cite{Hel}  it can be written as an alternating sum
\[  \int_{U_q} e^{-tr(\eta u \xi^2u^{-1}\eta)}du\,=\frac{\prod_{j=1}^{q-1}j!}{\pi(\xi^2)\pi(\eta^2)} \sum_{w\in S_q} sgn(w)\,e^{-(\xi^2, w\eta^2)}\]
where $(\,.\,,.\,)$ denotes the usual Euclidean scalar product in $\b R^q$ and 
$\,\pi(\xi) = \prod_{i<j} (\xi_i-\xi_j)\,$ is the fundamental alternating polynomial.

\end{enumerate}
\end{remarks}

\section{Laws of large numbers}

Let $\nu\in M^1(\Pi_q)$ be a probability measure and $\mu>\rho-1=d(q-1/2)$ 
a fixed index. We say that a time-homogeneous Markov chain $(S_k^\mu)_{k\ge0}$
on $\Pi_q$ is a Bessel-type random walk on $\Pi_q$ of index $\mu$ with law
$\nu$ if $S_0^\mu=0$ and if its transition probability is given by
$$P(S_{k+1}^\mu\in A| S_k^\mu=x) \> =\> (\delta_x*_\mu \nu)(A)$$
for all $k\in \mathbb N_0$, $x\in\Pi_q$ and Borel sets $A\subseteq \Pi_q$. 
It is easily checked by induction on $k$ that the distribution of
$S_k^\mu$
is just the $k$-fold convolution power $\nu^{(k,\mu)}=\nu *_{\mu}\nu *_{\mu}\ldots *_{\mu}\nu $ 
of $\nu$ with respect to the Bessel convolution of index $\mu$. As announced
in the introduction, we are interested in limit theorems for the random
variables $S_k^\mu$ as $k,\mu\to\infty$. Our first result in this direction  is the
following weak law of large numbers:

\begin{theorem}\label{LLN}
Let $\nu\in M^1(\Pi_q)$ be a probability measure with finite second moment 
$$\sigma^2(\nu):=\int_{\Pi_q} s^2\> d\nu(s) \in \Pi_q,$$
and let $(\mu_k)_{k\in\b N}\subset ]\rho-1,\infty  [$ be an arbitrary sequence
of indices with $\lim_{k\to\infty} \mu_k=\infty$. Let $S_k^{\mu_k}$ be the
$k$-th member of the  Bessel-type random walk of index $\mu_k$ with law $\nu$.
Then
$$\frac{1}{\sqrt k} S_k^{\mu_k}\to \sqrt{\sigma^2(\nu)}$$
in probability as $k\to\infty$.
\end{theorem}

This first main result  has the following consequence which was
stated as Theorem \ref{LLN2} in the introduction:

\begin{corollary}
Let $\nu\in M^1(\Pi_q)$ be a probability measure with finite second moment
$\sigma^2(\nu)\in \Pi_q$. For each dimension $p\in\b N$ consider the unique
$U_p$-invariant probability measure $\nu_p\in M^1(M_{p,q})$ with
$\phi_p(\nu_p)=\nu$ where $ \phi_p: M_{p,q} \to \Pi_{q},\> x\mapsto
(x^*x)^{1/2},$ is the 
canonical projection. Let further
$(X_l^p)_{l\in\b N}$ be a sequence of i.i.d.~$M_{p,q}$-valued random variables
with law $\nu_p$.
Then for each sequence $(p_k)_{k\in\b N}\subset \b N$ of dimensions with
$\lim_{k\to\infty} p_k=\infty$, the  $\Pi_{q}$-valued random variables 
$$\frac{1}{\sqrt k} \phi_{p_k}\bigl(\sum_{l=1}^k X_l^{p_k}\bigr)$$
tend in probability to the constant $\,\sqrt{\sigma^2(\nu)}$.
\end{corollary}

\begin{proof}
This is clear from Theorem \ref{LLN} because
$\bigl(\phi_{p}\bigl(\sum_{l=1}^k X_l^{p}\bigr)\bigr)_{k\ge0}$
is a Bessel-type random walk on $\Pi_q$ with index $\mu=pd/2$.
\end{proof}

The proof of Theorem  \ref{LLN} relies on estimates for matrix Bessel functions from the preceding
section and on   standard properties of the Laplace transform on matrix
cones. These properties are likely to be known but we include them for the reader's convenience.

Recall that the Laplace transform $L\nu\in C_b(\Pi_q)$ of a
measure $\nu\in M^1(\Pi_q)$ is defined by 
\[L\nu(x)=\int_{\Pi_q} e^{-\langle
  x,y\rangle} d\nu(y), \quad x\in \Pi_q.\]
The Laplace transform on the cone $\Pi_q$ satisfies the following
Levy-type continuity theorem.

\begin{proposition}\label{levy} For probability measures $\nu,\, (\nu_k)_{k\geq 1} \in  M^1(\Pi_q)$ the
  following statements are equivalent:
\begin{enumerate}\itemsep=-1pt
\item[\rm{(1)}] $\nu_k\to \nu$ weakly.
\item[\rm{(2)}] $L\nu_k(x)\to L\nu(x)$ for all $x\in\Pi_q$.
\item[\rm{(3)}] $L\nu_k(x)\to L\nu(x)$ for all $x\in\Omega_q$.
\end{enumerate}
\end{proposition}

\begin{proof}
$(1)\Longrightarrow (2)\Longrightarrow (3)$ is obvious. For $
(3)\Longrightarrow  (1)$ observe that for $x\in \Omega_q$, the
exponential function $e_x(y):=e^{-\langle  x,y\rangle}$ is contained in
$C_0(\Pi_q)$, i.e. it vanishes at infinity. Moreover, the linear span of $\{e_x, \,x\in \Omega_q\}$ is a
$\|.\|_\infty$-dense subspace of $C_0(\Pi_q)$ by the 
Stone-Weierstrass theorem.
It follows from (3) and a $3\epsilon$-argument that $\int f\> d\nu_k\to \int
f\> d\nu$ for all $f\in C_0(\Pi_q)$ which implies (1).
\end{proof}

The following result can be readily derived from the dominated convergence
theorem as in the classical setting:

\begin{lemma}\label{Laplace-Taylor1} Let $\nu\in  M^1(\Pi_q)$ be a probability
  measure which admits $r$-th moments for $r\in \b N$, i.e., 
$\int_{\Pi_q} \|y\|^r\> d\nu(r)<\infty$. Then $L\nu$ is $r$-times continuously
differentiable on $\Pi_q$.
\end{lemma}

Using the  Taylor formula at $0\in\Pi_q$, we in particular obtain:

\begin{corollary}\label{Laplace-Taylor2}
Let $\nu\in M^1(\Pi_q)$ with finite second moment $\sigma^2(\nu)\in \Pi_q$.
Then
\begin{equation}\label{Taylor}
\int_{\Pi_q} e^{-\langle sx, sx\rangle } \> d\nu(s) = 1- tr(x\sigma^2(\nu)x)
+o(\|x\|^2)\quad\text{ as }\, x\to 0\,\text{ in }\, \Pi_q.
\end{equation}
Moreover, if $\nu\in M^1(\Pi_q)$ admits fourth moments, then even $O(\|x\|^4)$ is true instead of $o(\|x\|^2)$ in  relation (\ref{Taylor}).
\end{corollary}

The following result is a variant of the preceding corollary:

\begin{lemma} \label{lemma4.6}
Let $\nu\in M^1(\Pi_q)$ with finite second moment $\sigma^2(\nu)\in \Pi_q$, and
let $(\mu_k)_{k\ge 1}\subset ]0,\infty[$ be as in Theorem \ref{LLN}. Then for each
$x\in \Pi_q$,
$$ \int_{\Pi_q} J_{\mu_k}\bigl( \frac{\mu_k}{k} xs^2x\bigr)\> d\nu(s)= 1-
\frac{1}{k}tr\bigl(x\sigma^2(\nu)x\bigr)\, + \,o(\|x\|^2/k) \quad \text{ as }\,k\to\infty. $$
Moreover,  if $\nu\in M^1(\Pi_q)$ admits fourth moments, then the error term
 $o(1/k)$ can be replaced by  $$O\bigl(\|x\|^4/k^{2} + \|x\|^2/(k\mu_k)\bigl).$$
\end{lemma}

\begin{proof} We first conclude from Theorem \ref{theorem1} that for $y\in\Pi_q$,
$$ \bigl|  J_{\mu_k}(\mu_k y)- e^{-tr\, y} \bigr| \le \frac{c}{\mu_k}\,tr\, y.$$
Therefore
\begin{align}
\int_{\Pi_q} \Bigl|  J_{\mu_k}\Bigl(\frac{\mu_k}{k} xs^2x) -
e^{-\frac{1}{k}tr(xs^2x)} \Bigr|\> d\nu(s) &\le
\frac{c}{k\mu_k}\int_{\Pi_q} tr(xs^2x)\> d\nu(s)
\notag\\
&\le \frac{c\|x\|^2}{k\mu_k}\int_{\Pi_q}\|s\|^2\> d\nu(s)
\notag\\
&\le \frac{\tilde c\|x\|^2}{k\mu_k}
\notag
\end{align}
 with suitable constants $c,\tilde c>0$. On the other hand, we conclude from
 Corollary \ref{Laplace-Taylor2} that
$$\int_{\Pi_q} e^{-\frac{1}{k}tr(xs^2x)} \> d\nu(s)=1- \frac{1}{k}tr(x\sigma^2(\nu)x)
+o(\frac{\|x\|^2}{k})$$
which yields the first claim.
The second statement follows readily from the second statement in  Corollary \ref{Laplace-Taylor2}.
\end{proof}

\begin{proof}[Proof of Theorem \ref{LLN}]  
Let $\nu^{(k,\mu_k)}$ be the
$k$-fold Bessel convolution power of $\nu$ with index $\mu_k$. Then
$\nu^{(k,\mu_k)}$ is the distribution of the random variable $S_k^{\mu_k}$.
Being hypergroup characters, the matrix Bessel functions $s\mapsto J_{\mu}(xs^2x)$ are
multiplicative w.r.t.~the Bessel convolution of index $\mu$. Together with the  preceding lemma this implies that
\begin{align}
\lim_{k\to\infty} \int_{\Pi_q} & J_{\mu_k}\Bigl(\frac{\mu_k}{k} xs^2x\Bigr)\> d\nu^{(k,\mu_k)}(s)
\notag\\
&= \lim_{k\to\infty} \Bigl( \int_{\Pi_q}  J_{\mu_k}\Bigl(\frac{\mu_k}{k}
xs^2x\Bigr)\> d\nu(s)\Bigr)^k
\notag\\
&= \lim_{k\to\infty} \Bigl( 1-\frac{1}{k} tr(x\sigma^2(\nu)x) +o(1/k)\Bigr)^k
\notag\\
&= e^{-tr(x\sigma^2(\nu)x)} =:A(x).
\end{align}
We thus conclude from Proposition \ref{proposition1} that
$$\lim_{k\to\infty} \int_{\Pi_q} e^{-\frac{1}{k}tr(xs^2x)}\> d\nu^{(k,\mu_k)}(s)
= \lim_{k\to\infty} \int_{\Pi_q}  J_{\mu_k}\Bigl(\frac{\mu_k}{k} xs^2x)\> d\nu^{(k,\mu_k)}(s)
=A(x)$$
for $x\in\Pi_q$. From this we conclude (after a quadratic transformation of the argument) that the Laplace transforms of the distributions of
$(S_k^{\mu_k})^2/k$ tend  to the Laplace transform of the point
measure $\delta_{\sigma^2(\nu)}$ on $\Pi_q$ as $k\to\infty$. The theorem now follows from
Proposition \ref{levy}.
\end{proof}

We next turn to a strong law of large numbers which generalizes Theorems
\ref{SLLN1} and \ref{SLLN2}.

\begin{theorem}\label{SLLN}
Let $\nu\in M^1(\Pi_q)$ be a probability measure with compact support.
 Let 
 $(\mu_k)_{k\in\b N}\subset ]\rho-1,\infty  [$ be an arbitrary sequence
of indices and
$(n_k)_{k\in\b N}\subset\b N$ a sequence of time steps with the
following properties:
\begin{enumerate}\itemsep=-1pt
\item[\rm{(1)}] $\lim_{k\to\infty}{\mu_k}/{k^a}=\infty$ for all $a\in\b N;$
\item[\rm{(2)}] $\lim_{k\to\infty}{\mu_k}/(n_k^2 (\ln k)^2)=\infty;$
\item[\rm{(3)}] $\lim_{k\to\infty}n_k/(\ln k)^2 =\infty$.
\end{enumerate}
Let $S_{n_k}^{\mu_k}$ be the
$n_k$-th member of the  Bessel-type random walk of index $\mu_k$ with law $\nu$.
Then,
$$\frac{1}{\sqrt{n_k}} S_{n_k}^{\mu_k}\to \sqrt{\sigma^2(\nu)}$$
for $k\to\infty$ almost surely.
\end{theorem}

As for the WLLN in the beginning of this section, this theorem immediately 
implies Theorems \ref{SLLN1} and  \ref{SLLN2}. 

Recall that the dimension of $H_q$ as a real vector space is given by $n= q + \frac{d}{2}q(q-1)$. The proof of Theorem \ref{SLLN} relies on the following elementary
observation:

\begin{lemma}\label{reduction-product}
There exist matrices $b_1,\ldots, b_n\in \Pi_q$ such that for all
 $a\in \Pi_q$ and sequences $(a_k)_{k\in \b N}\subset\Pi_q$  we have $a_k\to a$ if and only if
$\langle b_j,a_k\rangle \to \langle b_j,a\rangle $ for all $j=1,\ldots,n$.
\end{lemma}

\begin{proof}
If $b_1,\ldots, b_n$ is any $\b R$-basis of the vector space $H_q$ of Hermitian
matrices with  dimension $n=q+q(q-1)d/2$, then obviously  $a_k\to a$
if and only if $\langle b_j,a_k\rangle \to \langle b_j,a\rangle $ for all
$j=1,\ldots,n$.
On the other hand, we can find  a basis consisting of elements from $\Pi_q$.
For instance, we may take the $q$ diagonal matrices of the form
$diag(0,\ldots,0,1,0,\ldots,0)$  together with the matrices of the form
\begin{equation}
I+\frac{1}{2}\bigl( le_{i,j} + l^*e_{j,i}\bigr) \in \Pi_q 
\end{equation}
for $1\le i< j\le q$ and the $l\in\b F$ with $|l|=1$ forming an $\b R$-basis 
of $\b F$ where the $e_{i,j}$ are the elementary matrices with 1 in the
$(i,j)$-coordinate and $0$ otherwise. Notice that these matrices are positive
definite by the Gershgorin criterion.
\end{proof}

\begin{proof}[Proof of Theorem \ref{SLLN}]  
Let $\mu_k$ and $n_k$ be given as in the theorem.
By Lemma \ref{reduction-product}, it suffices to prove that for each 
$c\in\Pi_q$, 
\begin{equation}\label{scalar-bed}
\frac{1}{ {n_k}} \langle c^2, (S_{n_k}^{\mu_k})^2\rangle \to  \langle c^2,
\sigma^2(\nu)\rangle
\quad\quad
\text{ almost surely.}
\end{equation}
For this we shall prove for each $\epsilon>0$ that
\begin{equation}\label{scalar-bed1}
P\Bigl(\frac{1}{ {n_k}} \langle c^2, (S_{n_k}^{\mu_k})^2\rangle \ge \langle
c^2,\sigma^2(\nu)\rangle +\epsilon\Bigr) = O(1/k^2)
\end{equation}
and
\begin{equation}\label{scalar-bed2}
P\Bigl(\frac{1}{{n_k} } \langle c^2, (S_{n_k}^{\mu_k})^2\rangle \le \langle
c^2,\sigma^2(\nu)\rangle -\epsilon\Bigr) = O(1/k^2).
\end{equation}
Relation (\ref{scalar-bed}) then follows immediately from the Borel-Cantelli lemma.

\smallskip

We first turn to the proof of relation (\ref{scalar-bed2}). Here we proceed as in the
beginning of the proof of Lemma \ref{lemma4.6} and conclude from Theorem
\ref{theorem1}
that
\begin{align}\label{n1}
E\Bigl( e^{-\frac{2 \ln k}{\epsilon \cdot n_k}\cdot tr(c^2
  (S_{n_k}^{\mu_k})^2)}\Bigr)
&=\int_{\Pi_q}  e^{-\frac{2 \ln k}{\epsilon \cdot n_k}\cdot tr(cs^2c)}\> d\nu^{(n_k,\mu_k)}(s)
\notag\\
&= \int_{\Pi_q}  J_{\mu_k}\Bigl(\frac{2\mu_k \ln k}{\epsilon \cdot n_k}\cdot
cs^2c\Bigr)\> d\nu^{(n_k,\mu_k)}(s) +O\Bigl(\frac{1}{\mu_k}\Bigr).
\notag\\
&= \Bigl(\int_{\Pi_q}  J_{\mu_k}\Bigl(\frac{2\mu_k \ln k}{\epsilon \cdot n_k}\cdot
cs^2c\Bigr)\> d\nu(s) \Bigr)^{n_k} +O\Bigl(\frac{1}{\mu_k }\Bigr).
\end{align}
Moreover, using the stronger statement of Lemma \ref{lemma4.6},
Eq.~(\ref{mit0}), and the  assumptions (1) and (3) of the theorem, we obtain
\begin{align}\label{n2}
\Bigl(\int_{\Pi_q} & J_{\mu_k}\bigl(\frac{2\mu_k \ln k}{\epsilon \cdot n_k}\cdot
cs^2c\bigr)\> d\nu(s)\Bigr)^{n_k}
\notag\\
&= 
\Bigl(1-\frac{2 \ln k}{\epsilon \cdot n_k}\cdot tr(c\sigma^2(\nu)c) 
+O\Bigl(\frac{(\ln k)^2}{n_k^2}+ \frac{\ln k}{n_k\mu_k}\Bigr)\Bigr)^{n_k}
\notag\\
&\le e^{-\frac{2 \ln k}{\epsilon}  \cdot tr(c\sigma^2(\nu)c)} \cdot O(1)
\end{align}
The Markov inequality and estimates \eqref{n1}, \eqref{n2} now lead to
\begin{align}
P\Bigl(\frac{1}{n_k }& \langle c^2, (S_{n_k}^{\mu_k})^2\rangle \le \langle
c^2,\sigma^2(\nu)\rangle -\epsilon\Bigr)
\notag\\
&= P\Bigl( e^{-\frac{2 \ln k}{\epsilon \cdot n_k}\cdot tr(c^2
  (S_{n_k}^{\mu_k})^2)} \ge e^{-\frac{2 \ln k}{\epsilon }(
  tr(c\sigma^2(\nu)c)-\epsilon)}\Bigr)
\notag\\
&\le \frac{1}{ e^{-\frac{2 \ln k}{\epsilon }(
  tr(c\sigma^2(\nu)c)-\epsilon)}} 
\cdot E\bigl( e^{-\frac{2 \ln k}{\epsilon \cdot n_k}\cdot tr(c^2  (S_{n_k}^{\mu_k})^2)}\bigr)
\notag\\
&\le  e^{\frac{2 \ln k}{\epsilon } tr(c\sigma^2(\nu)c)} e^{-2\ln k} \Bigl(  e^{-\frac{2 \ln k}{\epsilon } tr(c\sigma^2(\nu)c)}
\cdot O(1) + O\Bigl(\frac{1}{\mu_k}\Bigr)\Bigr)
\notag\\
&\le O\Bigl(\frac{1}{k^2}\Bigr) + O\Bigl(\frac{k^a}{\mu_k}\Bigr)
\notag
\end{align}
with a suitable constant $a=a(\epsilon,c)>0$. Condition (1) of the theorem now completes the proof of
 (\ref{scalar-bed2}).

We now turn to the proof of relation (\ref{scalar-bed1}). Assume that
$supp\>\nu\subset \{x\in \Pi_q:\> \|x\|_2\le M\}$ holds for a suitable constant
$M>0$. 
Then by the support properties  of the Bessel convolution on $\Pi_q$, we have for all
$k\in\b N$
\begin{equation}\label{supportcond}
supp\>\nu^{(n_k,\mu_k)}\subset \{x\in \Pi_q:\> \|x\|_2\le n_kM\}.
\end{equation}
We now consider the function $H$ and the constant $C>0$ of Proposition
\ref{proposition3}.
We  conclude from Eq.~(\ref{supportcond}) and condition (2) of the theorem
that for all sequences $s_k\in supp\>\nu^{(n_k,\mu_k)}$,
$$\frac{(\ln k)^2\|s_k\|_2^4}{\mu_kn_k^2}\to 0 
\quad\quad\text{and}\quad 
  \frac{1}{\sqrt{\mu_k-\rho}} \cdot \sqrt{\frac{\ln k}{ n_k}}\cdot s_k\,\to 0.$$
Thus by the definition of $H$ we have for each $c\in \Pi_q$
$$H\Bigl(\sqrt{\frac{2\ln k}{\epsilon n_k}}\cdot cs_k, \sqrt{\mu_k-\rho}\Bigr)\to 0$$
and 
$$R_k(s):=\Bigl( 1- \frac{4C}{\epsilon^2}\cdot \frac{(\ln k)^2  \|s\|_2^4}{\mu_k n_k^2} - 
H\Bigl(\sqrt{\frac{2\ln k}{\epsilon n_k}}\cdot cs_k,
\sqrt{\mu_k-\rho}\Bigr)\Bigr)^{-1}$$
remains bounded as $k\to \infty$ and  $s\in supp\>\nu^{(n_k,\mu_k)}$.
This fact together with the estimates of Proposition \ref{proposition3} and
conditions (2) and (3) of the theorem imply that
\begin{align}\label{mainest2}
E\Bigl( e^{\frac{2 \ln k}{\epsilon \cdot n_k}\cdot tr(c
  (S_{n_k}^{\mu_k})^2c)}\Bigr)
&=\int_{\Pi_q}  e^{\frac{2 \ln k}{\epsilon \cdot n_k}\cdot tr(cs^2c)}\> d\nu^{(n_k,\mu_k)}(s)
\notag\\
&\le \int_{\Pi_q}  J_{\mu_k}\Bigl(-(\mu_k-\rho)\frac{2\ln k}{\epsilon \cdot n_k}\cdot
cs^2c\Bigr)\> d\nu^{(n_k,\mu_k)}(s) \cdot O(1)
\notag\\
&= \Bigl(\int_{\Pi_q}  J_{\mu_k}\Bigl((\mu_k-\rho)\frac{2\ln k}{\epsilon \cdot n_k}\cdot
cs^2c\Bigr)
\> d\nu(s) \Bigr)^{n_k}  \cdot O(1)
\notag\\
&\le  \Bigl(\int_{\Pi_q}  e^{\frac{2 \ln k}{\epsilon \cdot n_k}\cdot
  tr(cs^2c)}\> d\nu(s) \Bigr)^{n_k}
\cdot \bigl(1+C/\mu_k\bigr)^{n_k} \cdot O(1)
\notag\\
&\le  \Bigl(\int_{\Pi_q}  \Bigl(1+\frac{2 \ln k}{\epsilon \cdot n_k}\cdot
  tr(cs^2c)     +   O((\ln k/n_k)^2) \Bigr)\> d\nu(s) \Bigr)^{n_k}
  e^{Cn_k/\mu_k} \cdot O(1)
\notag\\
&= \Bigl(1+\frac{2 \ln k}{\epsilon \cdot n_k}\cdot tr(c\sigma^2(\nu)c)+ O((\ln k/n_k)^2) 
 \Bigr)^{n_k}
   \cdot O(1) 
\notag\\
&\le e^{\frac{2 \ln k}{\epsilon }\cdot tr(c\sigma^2(\nu)c)}  \cdot O(1).
\end{align}
Employing again the Markov inequality we thus obtain
\begin{align}
P\Bigl(\frac{1}{n_k }& \langle c^2, (S_{n_k}^{\mu_k})^2\rangle \ge \langle
c^2,\sigma^2(\nu)\rangle +\epsilon\Bigr)
\notag\\
&= P\Bigl( e^{\frac{2 \ln k}{\epsilon \cdot n_k}\cdot tr(c
  (S_{n_k}^{\mu_k})^2c)} \ge e^{\frac{2 \ln k}{\epsilon }(
  tr(c\sigma^2(\nu)c)+\epsilon)}\Bigr)
\notag\\
&\le \frac{1}{ e^{\frac{2 \ln k}{\epsilon }(
  tr(c\sigma^2(\nu)c)+\epsilon)}} 
\cdot E\bigl( e^{\frac{2 \ln k}{\epsilon \cdot n_k}\cdot tr(c  (S_{n_k}^{\mu_k})^2c)}\bigr)
\notag\\
&\le e^{-2\ln k}  \cdot O(1) \quad = \,O(1/k^2)
\notag
\end{align}
as claimed. This  proves Eq.~(\ref{scalar-bed1}) and completes the proof of the theorem.
\end{proof}

\begin{remarks}
\begin{enumerate}\itemsep=-1pt
\item[\rm{(1)}] Let us briefly comment on  the conditions of Theorem \ref{SLLN}. The
  most interesting case appears for $n_k=k$, where only the growth
  condition (1) on the indices $\mu_k$ and the compact support condition for $\nu$
  remain. Condition (1) is the essential condition in the end of the proof of
  Eq.~(\ref{scalar-bed2}), and we see no possibility to weaken this one.
On the other hand, the compact support of $\nu$ has been used mainly in order to
derive estimate (\ref{mainest2}) in a smooth way. We expect that here 
somewhat more involved estimations (for example, by using H\"olders inequality
in between) might
also lead to (\ref{mainest2}) under weaker conditions on the support of $\nu$.  It is however clear that any proof along our
approach will need that square-exponential moments of $\nu$ exist, i.e. $\int_{\Pi_q}
e^{tr(cs^2c)}\> d\nu(s)<\infty$ for all $c\in\Pi_q$. 
\item[\rm{(2)}] We expect that there exist also central limit theorems
  associated with the laws of large numbers above. In particular, the
  convergence of $\chi^2$-distributions to normal 
distributions for $q=1$ and convergence of Wishart distributions to multidimensional normal distributions for $q\geq 2$ 
  suggest that in a CLT normal distributions  appear as limits after taking squares
  after suitable renormalizations. 
\item[\rm{(3)}] Let us briefly return to the  case  $q=1$
  discussed in Theorems  \ref{LLN1} and \ref{SLLN1}. In this context one might ask for
  limit theorems for series of random walks on series of two-point homogeneous
  spaces where the number of steps and the dimensions of these spaces tend to
  infinity. For  spheres and projective spaces over  $\b F= \b R, \b C$ 
or  $\b H$, central limit theorems were given in \cite{V11} and references cited therein.
It should also  be interesting to study the non-compact cases, i.e. random walks on hyperbolic
spaces.
\item[\rm{(4)}] As explained  in Section 4, there is a
  close connection between Bessel convolutions on the matrix cones $\Pi_q$ and
  the theory of Dunkl operators on a $B_q$-Weyl chamber in $\b R^q$ for certain
  indices. It is clear that we may project Theorems  \ref{LLN} and  \ref{SLLN}
  to these particular cases.  We do not state this result separately.
Under the hypothesis that Dunkl operators are
  related to commutative hypergroups on Weyl chambers for all root systems and
  all positive multiplicities (see \cite{R3}), it will become an interesting question in Dunkl theory 
whether there exist laws of large numbers for random walks on Weyl chambers similar
  to Theorems \ref{LLN} and  \ref{SLLN} when the multiplicities of Dunkl theory tend to
  infinity.
\end{enumerate}
\end{remarks}

\section{A large deviation principle}

In this section we derive a large deviation principle (LDP) for $q=1$ and $\b
F =\b R$ which fits
to the laws of large numbers given in  Theorems  \ref{LLN1} and  \ref{SLLN1}. Before going into
details we explain the restriction  $q=1$.
Our proof of a LDP will  be based on  the limits
\begin{equation}
E\Bigl( e^{ \langle c,  (S_{n_k}^{\mu_k})^2\rangle}\Bigr)
\quad\quad\text{for}\quad
k\to\infty \quad\text{and all}\quad c\in H_{q}
\end{equation}
(in the notion of the preceding section) together with a standard result from
LDP theory (see e.g.~Theorem II.6.1 of Ellis \cite{E}) which states that  suitable
convergence of Laplace transforms implies a LDP. Unfortunately we can prove this
convergence only for matrices of the form $\pm c\in H_{q}$ with $c\in \Pi_q$, as our convergence
proofs
 depend  on estimates  for the Bessel
functions $J_\mu$ which were derived in Section 3 from the integral
representation (\ref{int-rep})  which is not available for
arbitrary matrices $c\in H_q$ for $q\ge 2$. We therefore restrict our attention to $q=1$
and consider the Bessel-type random walks    $(S_k^\mu)_{k\ge0}$
on $[0,\infty[=\Pi_1$ of indices $\mu$ with fixed law
$\nu\in M^1([0,\infty[)$.

\begin{proposition}\label{LDP-estimate} 
Let $\nu\in M^1(\Pi_q)$ be a probability measure with compact support.
 Let 
 $(\mu_k)_{k\in\b N}\subset ]1/2,\infty  [$ be a sequence
of indices and  $(n_k)_{k\in\b N}\subset \b N$ a sequence of time steps with
 $n_k\to\infty$
and
$\lim_{k\to\infty} e^{ an_k}/\mu_k =0$ for all $a>0$.
Then, 
$c_k(t) := \frac{1}{n_k} \ln  E(e^{t (S_{n_k}^{\mu_k})^2}) $
converges for $t\in\b R$ and $k\to\infty$ to
$$c(t):=\ln \Bigl(\int_0^\infty e^{t s^2}\> d\nu(s)\Bigr).$$
\end{proposition} 

\begin{proof}
We proceed  as the proof of Theorem  \ref{SLLN}. The case $t=0$ is trivial.
Now let $t>0$ and put $h(t):= \int_0^\infty e^{t s^2}\>d\nu(s)$.
Using Theorem \ref{theorem1} twice, 
we obtain  that
\begin{align}\label{exp-conver}
c_k(-t)&= \frac{1}{n_k} \ln\Bigl(\int_0^\infty e^{-t s^2}\>
d\nu^{(n_k,\mu_k)}(s) \Bigr)
\notag \\
&= \frac{1}{n_k} \ln\Bigl(\int_0^\infty J_{\mu_k}(\mu_k t s^2) \>
d\nu^{(n_k,\mu_k)}(s)+O(1/\mu_k) \Bigr)
\notag \\
&= \frac{1}{n_k} \ln\Bigl(\Bigl(\int_0^\infty J_{\mu_k}(\mu_k t s^2) \>
d\nu(s) \Bigr)^{n_k} +O(1/\mu_k) \Bigr)
\notag \\
&= \frac{1}{n_k} \ln\Bigl( (h(-t)+O(1/\mu_k))^{n_k} +O(1/\mu_k) \Bigr)
\notag\\
&= \ln\Bigl(h(-t)+ O(1/\mu_k)\Bigr) +  \frac{1}{n_k}\ln\Bigl(1+
\frac{1}{(h(-t)+O(1/\mu_k))^{n_k}\mu_k}\Bigr)
\quad\to\quad  c(-t)
\end{align}
by the convergence conditions of the theorem.
 Furthermore, we obtain from the estimations in Proposition \ref{proposition3}
 and with the
notions there that
\begin{align}
\int_0^\infty e^{t s^2}\>d\nu^{(n_k,\mu_k)}(s) 
&\ge (1+O(1/\mu_k))^{-1} \cdot \int_0^\infty J_{\mu_k}(-(\mu_k-3/2) t s^2) \>d\nu^{(n_k,\mu_k)}(s) 
\notag \\
&= (1+O(1/\mu_k))^{-1} \cdot\Bigl( \int_0^\infty J_{\mu_k}(-(\mu_k-3/2) t
s^2) \>d\nu(s)\Bigr)^{n_k}
\notag \\
&\ge \frac{1}{1+O(1/\mu_k)} \Bigl( \int_0^\infty e^{t s^2}\Bigl[1- Cs^4
t^2/\mu_k - H(s\sqrt t, \sqrt{\mu_k-3/2})\Bigr] \>d\nu(s)\Bigr)^{n_k}.
\notag
\end{align}
As $[\ldots]\to 1$ uniformly on the compact set $supp\>\nu$, it follows
readily that $\liminf c_k(t)\ge c(t)$. Finally, Proposition
\ref{proposition3}, $supp\>\nu^{(n_k,\mu_k)} \subset [0,Mn_k]$ for a suitable
$M>0$, and the convergence condition of the theorem imply
\begin{align}
\int_0^\infty e^{t s^2}\>d\nu^{(n_k,\mu_k)}(s) 
&\le  \int_0^\infty \frac{J_{\mu_k}(-(\mu_k-3/2) t s^2)}{1- Cs^4
t^2/\mu_k - H(s\sqrt t, \sqrt{\mu_k-3/2})} \>d\nu^{(n_k,\mu_k)}(s) 
\notag \\
&\le (1+o(1/n_k))  \int_0^\infty J_{\mu_k}(-(\mu_k-3/2) t s^2) \>d\nu^{(n_k,\mu_k)}(s) 
\notag \\
 &=(1+o(1/n_k))\Bigl( \int_0^\infty J_{\mu_k}(-(\mu_k-3/2) t s^2) \>d\nu\Bigr)^{n_k}
\notag \\
 &\le   (1+o(1/n_k))(1+O(1/\mu_k))^{n_k}  h(t)^{n_k}
\notag
\end{align}
and thus $\limsup c_k(t)\le c(t)$. In summary,   $c_k(t)\to c(t)$ for
$t>0$ which completes the proof. 
\end{proof}

Notice that the very strong convergence condition in the proposition was
needed in the end of (\ref{exp-conver}) only where $h(-t)<1$ may become arbitrarily
small. This is caused by the fact that the difference estimation in Theorem
\ref{theorem1} does not fit well to the "multiplicative" structure of LDPs.
For all other estimates in the proof above much weaker polynomial convergence conditions are sufficient.

We here notice that the free energy function $c$ of Proposition
\ref{LDP-estimate} is precisely the same as for the classical LDP of Cramer
for sums of i.i.d. random variables on $[0,\infty[$ with common law $\nu$ (see e.g.~Ch.~II.4
of   \cite{E}). Moreover,
Proposition \ref{LDP-estimate} together with Theorem II.6.1 of Ellis \cite{E}
immediately imply that in the setting of Proposition \ref{LDP-estimate}, the distributions  of the
random variables $(S_{n_k}^{\mu_k})^2$  have the large deviation property with
scaling parameters $n_k$ and the rate function
$$I(s):=\sup_{t\in \b R} (st-c(t))\quad\quad (s\in\b R)$$
in the sense of Definition II.3.1 of  \cite{E}. We skip the details here.

\begin{remark}
If the conditions of   Proposition
\ref{LDP-estimate} are satisfied, we obtain that the  free energy function $c$
is differentiable on $\b R$ with $c^\prime (0)= \sigma^2(\nu)>0$.  Theorems
II.6.3 and II.6.4 of  \cite{E} now imply that (after taking square roots)  $S_{n_k}^{\mu_k}/\sqrt{n_k}$ converges to
$\sqrt{\sigma^2(\nu)}$ almost surely. Notice that this strong law of large numbers (SLLN) holds
under conditions which are slightly different from those in Theorem \ref{SLLN}
for $q=1$. This SLLN can be also derived directly for arbitrary $q\ge1$ similar
to the proof of Theorem \ref{SLLN}. As the conditions concerning the
parameters $\mu_k$ are extremely strong here, we omit details.
\end{remark}

\end{document}